\documentclass[12pt]{elsart}
\usepackage{amsfonts,amssymb}
\journal{Advances in Mathematics}

\newtheorem{Proposition}{Proposition}
\newtheorem{Theorem}{Theorem}
\newtheorem{Lemma}{Lemma}

\newenvironment{Proof}{{\bf Proof:}}{\hfill\rule{2mm}{2mm}\abst}

\newcommand{\abst}{\vspace{2ex}}

\begin{document}
\begin{frontmatter}
\title{Level Crossing Probabilities I: One-dimensional Random Walks and
Symmetrization}
\author{Rainer Siegmund-Schultze and Heinrich von Weizs\"{a}cker}
\address{Technische Universit\"{a}t Kaiserslautern\\
Fachbereich Mathematik\\
Erwin-Schr\"{o}dinger-Str., Geb\"{a}ude 48\\
67663 Kaiserslautern, Germany}

\begin{abstract}
We prove for an arbitrary random walk in $\mathbf{R}^1$ with independent
increments that the probability of crossing a level at a given time $n$ is $%
O(n^{-1/2})$. Moment or symmetry assumptions are not necessary. In removing
symmetry the (sharp) inequality $P(|X+Y| \le 1) < 2 P(|X-Y| \le 1)$ for
independent identically distributed $X,Y$ is used. In part II we shall
discuss the connection of this result to 'polygonal recurrence' of
higher-dimensional walks and some conjectures on directionally reinforced random walks
in the sense of Mauldin, Monticino and v.Weizs\"{a}cker \cite
{Mauldin-Monticino-Weizsaecker}.
\end{abstract}

\end{frontmatter}

\section{Introduction}

For a one-dimensional random
walk with independent steps it is a classical question to estimate the
probabilities of hitting a given level, to visit a given site or to decide
between recurrence and transience. In this context it is worthwhile to ask
for the probability that, for a given level $c$, this level is crossed (from
above or below) at a given time $n$. Denoting by $S_n=\sum_{i=1}^nX_i$ the
position at time $n$, with the increments $X_i$ being i.i.d., we say that
the random walk crosses the level $l$ at $n$, if sgn$(S_n-l)\neq $sgn$%
(S_{n-1}-l)$.

It is well-known that for a mean-zero random walk in $\mathbf{Z}^1$ with a
finite second moment the probability of hitting zero is exactly of the
order $n^{-1/2}$ if the walk is strongly aperiodic (\cite{Spitzer}, pp.42, 79), and is $O(n^{-1/2})$ in the general situation (except for
the trivial non-walk with $X_i=0$ a.s.) (\cite{Spitzer}, p.72). In the
general situation level crossings in the sense defined above can be much
more likely. For instance, in the case of a transient symmetric random walk,
the sum $\sum_{i=1}^\infty P(S_i=0)$ is finite, whereas by the fact that any
non-trivial symmetric random walk performs infinitely many changes of sign,
the Borel-Cantelli lemma shows that $\sum_{i=1}^\infty P($sgn$(S_i)\neq $sgn$%
(S_{i-1}))$ is a.s. infinite. So one might conjecture that for symmetric
random walks the probability of changing the sign at time $n$ could be much
larger than $n^{-1/2}.$

But this turns out to be wrong: Without any moment assumptions it can be
shown that level crossings are subject to the same $O(n^{-1/2})$-bound as
hitting probabilities (section 4, Theorem \ref{theo}).

The proof of this fact is surprisingly short in the symmetric case (section
2), where also a lower bound of the order $n^{-1}$ can be derived. Both
estimates make use of a few combinatorial arguments which essentially was introduced already by Erd\"{o}s and Hunt in \cite{Erdoes-Hunt} to
derive bounds in the symmetric case of the same order as those given here
(they additionally assume continuity of the random variables). The elegant
formulation of the combinatorial argument leading to the lower $\frac{1}{n}$
bound quoted here is due to Y. Peres.

In the general non-symmetric situation we make use of a new estimate
coupling an arbitrary increment distribution with its symmetrization
(section 3, Theorem \ref{symineqth}).

\section{Sign changes for symmetric random walks}

Let $\left( X_i\right) _{i=1,2,...}$ be an i.i.d. sequence of arbitrary
random variables in $\mathbf{R}^1$ with a symmetric distribution, i.e. $\mathcal{L}(X_i)= 
\mathcal{L}(-X_i).$ Consider the one-dimensional random walk $%
S_n=\sum_{i=1}^nX_i.$ Then we have

\begin{Proposition}
\label{sym} For a symmetric random walk $(S_{n})$ the probability of
changing sign at time $n$ is $O(n^{-1/2}):$ 
\begin{equation}\label{On^halb}
P({\rm sgn}(S_{n})\neq {\rm sgn}(S_{n-1}))=O(n^{-1/2})
\end{equation}
and we have the lower estimate: 
\begin{equation}\label{lowest}
P({\rm sgn}(S_{n})\neq {\rm sgn}(S_{n-1}))\geq \frac{1}{2n}\left(
1-P(X_{1}=0)^{n}\right) .
\end{equation}
This is sharp in the sense that there is a symmetric distribution with 
\begin{equation}\label{sharp}
P({\rm sgn}(S_{n})\neq {\rm sgn}(S_{n-1}))=\frac{1}{2n}+o(\frac{1}{n}).
\end{equation}
\end{Proposition}

\begin{Proof}
1.We rewrite $S_{n}$ as $S_{n}=\sum \varepsilon _{i}Y_{i}$, where the $%
Y_{i}=|X_{i}|$ are i.i.d. non-negative, the $\varepsilon _{i}$ are $\pm 1$
with equal probability and independent of each other and of the $Y_{i}.$
Then we have 
\begin{eqnarray*}
P({\rm sgn}(S_{n})\neq {\rm sgn}(S_{n-1})) &\leq &P(|S_{n-1}|\leq Y_{n})\
=\ 1-2P(S_{n-1}>Y_{n}) \\
&=&1-2P\left( \sum\limits_{i=1}^{n-1}\varepsilon _{i}Y_{i}-Y_{n}>0\right)  \\
&=&1-2P\left( {\rm sgn}\left( \sum\limits_{i=1}^{n-1}(-\varepsilon
_{i}\varepsilon _{n})Y_{i}+\varepsilon _{n}Y_{n}\right) =-\varepsilon
_{n}\right)  \\
&=&1-2P\left( {\rm sgn}\left( S_{n}\right) =-\varepsilon _{n}\right) 
\end{eqnarray*}
since $-\varepsilon _{i}\varepsilon _{n}$ has, for $i=1,2,...,n-1,$ the same
distribution as $\varepsilon _{i}$, and all these values are independent. We
continue the chain of inequalities 
\begin{eqnarray*}
&&1-2P\left( {\rm sgn}\left( S_{n}\right) =-\varepsilon _{n}\right) =1-%
\frac{2}{n}\mathbf{E}\#\left\{ i:{\rm sgn}\left( S_{n}\right) =-\varepsilon
_{i}\right\}  \\
&=&2\mathbf{E}\left( \frac{1}{2}-\frac{1}{n}\cdot \#\{i:{\rm sgn}\left(
S_{n}\right) =-\varepsilon _{i}\}\right)  \\
&\leq &2\mathbf{E}\left| \frac{1}{2}-\frac{1}{n}\cdot \#\{i:{\rm sgn}\left(
S_{n}\right) =-\varepsilon _{i}\}\right|  \\
&\leq &2P\left( S_{n}=0\right) +2\mathbf{E}\left| \frac{1}{2}-\frac{1}{n}%
\cdot \#\{i:1=\varepsilon _{i}\}\right| +2\mathbf{E}\left| \frac{1}{2}-\frac{%
1}{n}\cdot \#\{i:-1=\varepsilon _{i}\}\right|  \\
&\leq &2P\left( S_{n}=0\right) +4\sqrt{\mathbf{E}\left( \frac{1}{2}-\frac{1}{%
n}\cdot \#\{i:1=\varepsilon _{i}\}\right) ^{2}} \\
&=&2P\left( S_{n}=0\right) +2\sqrt{\frac{1}{n}},
\end{eqnarray*}
since the expression in the
last but one line is simply the square root of the variance of a binomial
distribution. Now (except for the trivial case $X_{n}\equiv 0$) $P\left( S_{n}=0\right) $ is $O(n^{-\frac{1}{2}})$, see 
Theorem 3 of \cite{Esseen}. This proves the upper bound.

2. We have the following estimate 
\begin{eqnarray*}
&&P({\rm sgn}(S_{n})\neq {\rm sgn}(S_{n-1}))=P({\rm sgn}(S_{n})\neq {\rm %
sgn}(S_{n}-\varepsilon _{n}Y_{n})) \\
&=&\frac{1}{n}\sum\limits_{j=1}^{n}P({\rm sgn}(S_{n})\neq {\rm sgn}%
(S_{n}-\varepsilon _{j}Y_{j})) \\
&=&\frac{1}{2n}\sum\limits_{j=1}^{n}2\mathbf{E(1}_{\{{\rm sgn}(S_{n})\neq 
{\rm sgn}(S_{n}-\varepsilon _{j}Y_{j})\}}\mathbf{)} \\
&=&\frac{1}{2n}\mathbf{E}\sum\limits_{j=1}^{n}(\mathbf{1}_{\{{\rm sgn}%
(S_{n}^{(j)})\neq {\rm sgn}(S_{n}^{(j)}-\varepsilon _{j}Y_{j})\}}\mathbf{+%
\mathbf{1}}_{\{{\rm sgn}(S_{n}^{(j+1)})\neq {\rm sgn}(S_{n}^{(j+1)}+%
\varepsilon _{j}Y_{j})\}}\mathbf{)}
\end{eqnarray*}

Here we used the abbreviation $S_{n}^{(j)}$ for $\sum_{i=1}^{j-1}(-%
\varepsilon _{i})Y_{i}+\sum_{i=j}^{n}\varepsilon _{i}Y_{i}$ , $j=1,2,...,n+1$%
. Obviously $S_{n}^{(j)}$ has the same distribution as $S_{n},$ and we have $%
S_{n}^{(1)}=S_{n}=-S_{n}^{(n+1)}.$ If we assume for the moment that not all $%
Y_{i}$ are zero, then obviously there is some index $j_{0}$ , $1\leq
j_{0}\leq n$ with sgn$(S_{n}^{(j_{0})})\neq $sgn$(S_{n}^{(j_{0}+1)})$, and
we have $S_{n}^{(j_{0}+1)}=S_{n}^{(j_{0})}-2\varepsilon _{j_{0}}Y_{j_{0}}.$
So one of the quantities $S_{n}^{(j_{0})},S_{n}^{(j_{0}+1)}$ necessarily has
a different sign if compared with $S_{n}^{(j_{0})}-\varepsilon
_{j_{0}}Y_{j_{0}}=S_{n}^{(j_{0}+1)}+\varepsilon _{j_{0}}Y_{j_{0}}.$ So we
may conclude that 
\begin{eqnarray*}
P({\rm sgn}(S_{n})\neq {\rm sgn}(S_{n-1})) &\geq &\frac{1}{2n}\left(
1-P(Y_{i}=0,i=1,2,...,n)\right)  \\
&=&\frac{1}{2n}\left( 1-P(Y_{1}=0)^{n}\right) .
\end{eqnarray*}
This proves the lower bound.\newline
3. The lower bound is sharp: Let $p=(p_{k})_{k\in \mathbb{N}}$ be a
probability distribution on $\mathbb{N}$ such that $P(A_{n})=o(\frac{1}{n})$
where $A_{n}$ is the event that the largest and the second largest of $n$
independent samples of $p$ coincide. For example it is not difficult to
verify that this holds if $p_{k}={\rm const}\cdot k^{-3/2}$. Let $%
P(Y_{i}=k!)=p_{k}$ for all $k\in \mathbb{N}$. Then on $A_{n}^{c}$ the
maximal sample dominates the sum of all others. Hence 
\begin{eqnarray*}
P(\textrm{sgn}(S_{n})\not=\textrm{sgn}(S_{n-1})) &\leq &\frac{1}{2}P(Y_{n}\geq
|\sum_{j=1}^{n-1}\varepsilon _{j}Y_{j}|)+\frac{1}{2}P(S_{n-1}=0) \\
&\leq &\frac{1}{2}P(A_{n}^{c}\cap \{Y_{n}=\max_{j\leq n}Y_{j}\})+\frac{1}{2}%
P(A_{n})+o(\frac{1}{n}) \\
&=&\frac{1}{2n}P(A_{n}^{c})+o(\frac{1}{n})=\frac{1}{2n}+o(\frac{1}{n})
\end{eqnarray*}
where we have used that $P(S_{n-1}=0)\leq P(A_{n-1})=o(\frac{1}{n})$.
\end{Proof}

\textbf{Remark}. If in addition to symmetry the distribution function of the 
$X_{n}$ is continuous, the proof shows that the probability of sign change
is at most $2n^{-1/2}$; in particular the expectation of the number $N_{n}$
of sign changes up to time $n$ satisfies $\mathbb{E}(N_{n})\leq
2\sum_{k=1}^{n}k^{-1/2}$. P. Erd\"{o}s and G.A. Hunt (\cite{Erdoes-Hunt})
gave another upper estimate for the probability of a sign change and for $%
\mathbb{E}(N_{n})$ which implies that in this last estimate the constant 2
can be changed to $(8\pi )^{-1/2}+\varepsilon $ for large $n$. They also
gave a.s. results for the asymptotic behaviour of $N_{n}$.

\section{The symmetrization inequality}

We prove the following

\begin{Theorem}
\label{symineqth} \label{prop}Let $X,Y$ be independent real random variables
with the same distribution. Then for any $c>0$ we have 
\begin{equation}
P(|X+Y|\leq c)<2P(|X-Y|\leq c).  \label{symineq}
\end{equation}
\end{Theorem}

This inequality is optimal in the sense that for any $\gamma<2$ there exists
an example of a probability distribution such that $P(|X+Y|\le 1)>\gamma
\cdot P(|X-Y|\le 1).$ It is rather easy to prove the result for $\gamma=3$
and this would be sufficient for our purposes, but we prefer to demonstrate
the optimal inequality.

This section is influenced by Jochen Vo{\ss}, Tobias Wahl (both
Kaiserslautern) and later by Y. Peres. The first suggested to first study
counting measures on finite sets, the second gave a discrete analogue of the
proof of Lemma \ref{symineqlem} below. Y. Peres pointed out the similarity
to the 123 theorem and its generalizations in \cite{Alon-Yuster} and asked
whether strict inequality holds. The present form of Lemma \ref{symineqlem}
actually is very close in spirit to the 'claim' on p. 325 in \cite
{Alon-Yuster}. Our Proposition \ref{symineqprop} gives a general procedure
showing how to derive two-variable inequalities from one-variable estimates like the
lemma. In a proof of the 123 theorem one would apply it to the function $%
f(x,y) = 3\mathbf{1}_{\{ |x-y|\le 1\} } - \mathbf{1}_{\{ |x-y|\le 2\} } .$

Note that in higher dimensions these questions become more complicated and a
simple general characterization of those functions $f$ which satisfy the
second alternative in the proposition seems out of reach. A couple of 
interesting multidimensional versions of our Theorem (but without optimal 
constants) were 
given in \cite{Wiegand}. Finally it should
be pointed out that the mere existence of a symmetrization constant in higher
dimensions (replacing the 2 of theorem \ref{sym}) follows eg. from \cite
{Mattner}, estimate (29) applied to $\mu = \mathcal{L}(X), \nu = \mathcal{L}%
(-X)$, and $\psi(x) = \mathbf{1}_{\{\|x\| \le 1\}}$.

The optimality of the constant $2$ can be seen by the following example: Let 
$X,Y$ be equidistributed on $\{-2n+1,-2n+3,\cdots ,-1,2,4,\cdots ,2n\}$ and $%
c=1.5$. Then $P(|X-Y|\leq 1.5)=P(X=Y)=\frac{1}{2n}$ and $P(|X+Y|\leq
1.5|X=k)=P(Y\in \{-k-1,-k+1\})=\frac{1}{n}$ except for the cases $k=-1$ and $%
k=2n.$ So we have $P(|X+Y|\leq 1.5)\geq \frac{1}{n}(1-\frac{1}{n})=2\cdot 
\frac{1}{2n}+o(\frac{1}{n})$. Actually let $\gamma (\mu )$ denote the
infimum of all constants by which one can replace the 2 in (\ref{symineq}) if 
$X$ and $Y$ have the law $\mu $. A closer inspection of the proof of Lemma \ref
{symineqlem} below suggests that such a periodic pattern as in this example
necessarily appears asymptotically in the distributions $\mu _{n}$ whenever $%
\gamma (\mu _{n})$ approaches $2$.
\newline
\newline
\begin{Lemma}
\label{symineqlem} Let $P$ be a probability measure on $\mathbb{R}.$ Define
a function $p$ by $p(x)=P[x-1,x+1].$ Then 
\begin{equation}
P\{x:p(-x)<2p(x)\}>0  \label{positiv}
\end{equation}
\end{Lemma}

\begin{Proof}
Assume that (\ref{positiv}) is false, i.e. $P(A)=1$ for $A=\{x:p(-x)\geq
2p(x)\}.$ Let $\alpha =\sup_{x\in A}p(x)$, which is easily seen to be positive. In fact, represent $\mathbb{R}$ as a union of disjoint half-open intervals of length $1$. One of these, say $I$, has positive $P$-measure. Each $x\in I$ fulfils $p(x)\geq P(I)$, and since $P(A)=1$ we conclude $P(A\cap I)>0$, so $A\cap I$ is non-empty, and hence $\alpha \geq P(I)$.

Let $\varepsilon >0.$ We show that for every $x\in A$
with $p(x)\geq \alpha -\varepsilon $ there is some $x^{\prime} \in A$ with $x+1\leq x^{\prime}
\leq x+2$ and $P(x+1,x^{\prime} +1]\geq \alpha -4\varepsilon $.  Since $x\in A$ we
get $p(-x)\geq 2p(x)\geq 2\alpha -2\varepsilon .$ Since (in view of $P(A)=1$%
) 
\begin{equation}\nonumber
P([-x,-x+1])\leq \sup \{p(z):-x\leq z\leq -x+1,\ \,z\in A\}\leq \alpha 
\end{equation}
we conclude $P([-x-1,-x))\geq \alpha -2\varepsilon $ and hence there is some 
$y\in A$ with $-x-1 \leq y<-x$ and $p(y)\geq \alpha -2\varepsilon .$ Then $p(-y)\geq
2\alpha -4\varepsilon .$ On the other hand $x<-y\leq x+1$ and $p(x)\leq
\alpha $, i.e. $P(x+1,-y+1]\geq \alpha -4\varepsilon .$ In particular there
is some $x^{\prime }\in A$ with $x+1\leq x^{\prime } \leq -y+1\leq x+2$ and 
\begin{equation}\nonumber
p(x^{\prime })\geq P(x+1,\ x^{\prime }+1]\geq \alpha -4\varepsilon .
\end{equation}
Proceeding recursively we can construct a sequence $x=x_{0},x_{1},x_{2},%
\ldots $ such that $x_{n}+1\leq x_{n+1}\leq x_{n}+2$ and 
\begin{equation}\nonumber
p(x_{n+1})\geq P(x_{n}+1,x_{n+1}+1]\geq \alpha -4^{n+1}\varepsilon .
\end{equation}
The intervals $(x_{n}+1,x_{n+1}+1]$ are disjoint and hence 
\begin{equation}\nonumber
1\geq \sum_{k=0}^{n-1}P(x_{k}+1,x_{k+1}+1]\geq n\alpha -\varepsilon
\,\sum_{k=1}^{n}4^{k}.
\end{equation}

Choosing first $n$ large enough and then $\varepsilon $ small enough we
arrive at a contradiction.
\end{Proof}

\begin{Proposition}
\label{symineqprop} Let $(\Omega ,\mathcal{B})$ be a measurable space and
let $f:\Omega ^{2}\longrightarrow \mathbb{R}$ be a $\mathcal{B}\otimes 
\mathcal{B}$ measurable bounded symmetric function. Let $\mathcal{P}$ be the
set of all probability measures on $\mathcal{B}.$ Then the following
dichotomy holds: Either 
\begin{equation}\nonumber
\int_{\Omega }\,f(\cdot ,y)\ P(dy)\leq 0\ \ P-a.e.
\end{equation}
for some $P\in \mathcal{P}$ or 
\begin{equation}\nonumber
\int_{\Omega \times \Omega }f(x,y)\ P(dx)P(dy)>0
\end{equation}
for all $P\in \mathcal{P}.$
\end{Proposition}

\begin{Proof}
We consider the cone $\mathcal{M}=\mathbb{R}_{+}\mathcal{P}$ of positive
finite measures on $\mathcal{B}$ and the bounded symmetric bilinear form $Q:%
\mathcal{M}\times \mathcal{M}\longrightarrow \mathbb{R}$ given by 
\begin{equation}\nonumber
Q(\mu ,\nu )=\int \,f(x,y)\,\mu (dx)\nu (dy).
\end{equation}
Suppose the first alternative is not valid. Then for all $\mu \in \mathcal{M}
$ 
\begin{equation}
\mu \{x:\int \,f(x,y)\ \mu (dy)>0\}>0.  \label{positiv2}
\end{equation}
1. We prove $Q(\mu ,\mu )>0$ for all measures $\mu $ of the form $\mu
=\sum_{i=1}^{n}p_{x_i}\delta _{x_{i}}.$ In this case the bilinear form is given by
a symmetric real matrix $A=(a_{i,j})_{i,j=1}^{n}$ with $%
a_{i,j}=f(x_{i},x_{j}).$ We have to prove that if $Ap$ has a positive entry
at some index $i$ with $p_{i}>0$ for each (probability) vector $p\geq 0,\sum
p_{j}=1$, then necessarily $q^{\prime }Aq$ is positive for each probability
vector $q.$ In fact, consider the minimum of $q^{\prime }Aq$ over the
simplex $q\geq 0,\sum q_{j}=1.$ We may assume without any loss of generality that this minimum is attained at some $q^{\ast }>0$, since
otherwise we may simply pass to the quadratic submatrix $(a_{i,j})_{q_{i}^{%
\ast }>0,q_{j}^{\ast }>0}$ of $A$\ formed by all indices where $q^{\ast }$
is non-zero. Now the projection of the gradient $\nabla _{q}(q^{\prime
}Aq)=2Aq$ to the linear subspace given by the condition $\sum
q_{j}=(1,1,...,1)\cdot q=1$ must be zero at $q^{\ast }$, which means that $%
Aq^{\ast }$ is a multiple of $(1,1,...,1)^{\prime }.$ By assumption $%
Aq^{\ast }$ has a positive entry, so all entries are positive. This proves
that $q^{\ast \prime }Aq^{\ast }>0.$ 

2. We prove $Q(\mu ,\mu )\geq 0$ for every $\mu \in \mathcal{M}.$ Without
loss of generality we may assume that $\mu $ is a probability measure. Let $%
n\in \mathbb{N}$ and $X_{1},\ldots ,X_{n}$ be iid. with law $\mu .$ Then 
\begin{eqnarray*}
&&n\ \mathbb{E}_{\mu }(f(X,X))+n(n-1)Q(\mu ,\mu ) \\
&=&\sum_{i,j=1}^{n}\mathbb{E}(f(X_{i},X_{j}))=\mathbb{E}%
\sum_{i,j=1}^{n}f(X_{i},X_{j}) \\
&=&\mathbb{E}(Q(\sum_{i=1}^{n}\delta _{X_{i}},\sum_{i=1}^{n}\delta
_{X_{i}}))>0.
\end{eqnarray*}
Since this holds for all $n$ we must have $Q(\mu ,\mu )\geq 0.$

3. Finally we prove the strict inequality $Q(\mu ,\mu )>0.$ Assume the
opposite that $Q(\mu ,\mu )=0$ for some probability measure $\mu $. By
assumption we know that the set $B^{+}=\{x:\int f(x,y)\mu (dy)>0\}$ has
positive $\mu $-measure. This together with $Q(\mu ,\mu )=0$ would imply
that $B^{-}=\{x:\int f(x,y)\mu (dy)<0\}$ has positive $\mu $-measure, too.
Let $\nu :=\mu (\cdot |B^{-})$. We then have $Q(\mu ,\nu )<0$. Obviously $%
\mu _{t}:=(1-t)\mu +t\nu =\mu +t(\nu -\mu )$ is a probability measure for $%
0\leq t\leq 1$. Consider the real non-negative quadratic function  $\lambda
(t):=Q(\mu _{t},\mu _{t})=2tQ(\mu ,\nu )+t^{2}(Q(\nu ,\nu )-2Q(\nu ,\mu )).$
Since $\lambda (0)=0$ we conclude that $0\leq \lambda ^{\prime }(0)=2Q(\mu
,\nu )$, arriving at a contradiction.
\end{Proof}

\begin{Proof}
(of Theorem \ref{symineqth}). It is sufficient to consider the case $c=1$.
Then the result follows from Proposition \ref{symineqprop} applied to the
function $f(x,y)=2\mathbf{1}_{\{|x-y|\leq 1\}}-\mathbf{1}_{\{|x+y|\leq 1\}}$
since Lemma \ref{symineqlem} shows that the first alternative in the
proposition is not satisfied.
\end{Proof}

\section{The general level crossing estimate}

We are now in a position to prove the announced level crossing estimate for
arbitrary random walks.

\begin{Theorem}
\label{theo} For any one-dimensional random walk $(S_{n})$ the probability
of a crossing of level $l$ at time $n$ is $O(n^{-\frac{1}{2}})$: 
\begin{equation}
P({\rm sgn}(S_{n}-l)\neq {\rm sgn}(S_{n-1}-l))=O(n^{-\frac{1}{2}}).
\label{Ohalb}
\end{equation}
\end{Theorem}

\begin{Proof}
1. For the trivial random walk $X_{1}=0$ a.s. the assertion is trivial, too.
First we prove the result for non-trivial random walks with $X_{1}\geq 0$
a.s. (or $X_{1}\leq 0$ a.s.). In this case it is almost trivial: We have $P($%
sgn$(S_{n}-l)\neq $sgn$(S_{n-1}-l))\leq P(S_{n-1}\leq l)$, and this last
expression tends to zero even exponentially fast by standard large deviation
theory. So we may assume for the following that 
\begin{equation}
P(X_{1}>0)\cdot P(X_{1}<0)>0.  \label{+-}
\end{equation}

2. We show that it is sufficient to prove the result for $l=0$, i.e. for
sign changes. In fact, suppose we know that the probability of a sign change
is $O(n^{-\frac{1}{2}}).$ It is enough to consider the case $l>0$ . If there
is a sequence $n_{i}$ with $P($sgn$(S_{n_{i}}-l)\neq $sgn$%
(S_{n_{i}-1}-l))\cdot n_{i}^{\frac{1}{2}}\rightarrow \infty $, then by
assumption we may conclude that 
\begin{equation}\nonumber
P({\rm sgn}(S_{n_{i}}-l)\neq {\rm sgn}%
(S_{n_{i}-1}-l),S_{n_{i}}>0,S_{n_{i}-1}>0)\cdot n_{i}^{\frac{1}{2}%
}\rightarrow \infty 
\end{equation}
and hence $(P(S_{n_{i}}\in (0,l])+P(S_{n_{i}-1}\in (0,l]))\cdot n_{i}^{\frac{%
1}{2}}\rightarrow \infty $. This is in contradiction to the fact that
bounded sets have $O(n^{-\frac{1}{2}})$ probabilities, according to Theorem 3 of \cite{Esseen}. So in the sequel we
may confine ourselves to the case $l=0.$

3. We have 
\begin{equation}\nonumber
P\left( {\rm sgn}(S_{n})\neq {\rm sgn}(S_{n-1})\right) \leq P\left(
|S_{n-1}|\leq |X_{n}|\right) .
\end{equation}
So it is sufficient to show $P\left( |S_{n}|\leq |X|\right) =O(n^{-\frac{1}{2%
}})$, where $X$ has the same distribution as $X_{1}$ and is independent of $%
\{X_{1},...,X_{n}\}$. If we can prove this for even $n$, it follows for all $%
n$, since by (\ref{+-}) with a positive probability $X_{2n}$ has the
opposite to $S_{2n-1}$. In turn, to derive $P\left( |S_{2n}|\leq
|X|\right) =O(n^{-\frac{1}{2}})$, it is sufficient to show $P\left(
|S_{2n}|\leq |X-X^{\prime }|\right) =O(n^{-\frac{1}{2}})$ (where again $%
X^{\prime }$ has the same distribution as $X_{1}$ and is independent of $%
\{X_{1},...,X_{n},X\})$, since with a positive probability $X^{\prime }$ has
the opposite sign to $X$. Now observe that $S_{2n}$ can be written as $%
S_{n}+S_{n}^{\prime }$ with $S_{n}^{\prime }$ being an independent copy of $%
S_{n}$. Now by Theorem \ref{symineqth} we have $P\left( |S_{n}+S_{n}^{\prime
}|\leq |X-X^{\prime }|\right) \leq 2\cdot P\left( |S_{n}-S_{n}^{\prime
}|\leq |X-X^{\prime }|\right) $, and the last expression is simply $2\cdot
P\left( |S_{n}^{\ast }|\leq |X^{\ast }|\right) $ where $S_{n}^{\ast }$ is a
symmetric random walk with increments distributed as $X^{\ast }=X-X^{\prime }
$. This random walk is non-trivial by (\ref{+-}), and in the proof of
Proposition \ref{sym} it was shown that in that case $P\left( |S_{n}^{\ast
}|\leq |X^{\ast }|\right) $ is $O(n^{-\frac{1}{2}})$.
\end{Proof}

\textbf{Acknowledgement}. The authors are grateful to Yuval Peres for
several helpful comments and his continuing interest in this work. We wish
to thank Rongfeng Sun for pointing us to an error in the proof of 
Proposition 2 in a previous version of this article.

\end{document}